\documentclass[12pt]{amsproc}
\usepackage{amsmath,amssymb,amscd}

\newcommand{\cA}{{\mathcal A}}
\newcommand{\cC}{{\mathcal C}}
\newcommand{\cE}{{\mathcal E}}
\newcommand{\cF}{{\mathcal F}}
\newcommand{\cG}{{\mathcal G}}
\newcommand{\cH}{{\mathcal H}}
\newcommand{\cJ}{{\mathcal J}}
\newcommand{\cK}{{\mathcal K}}
\newcommand{\cL}{{\mathcal L}}
\newcommand{\cP}{{\mathcal P}}
\newcommand{\cR}{{\mathcal R}}
\newcommand{\bbC}{{\mathbb C}}
\newcommand{\bbN}{{\mathbb N}}
\newcommand{\bbR}{{\mathbb R}}
\newcommand{\bbT}{{\mathbb T}}
\newcommand{\bbZ}{{\mathbb Z}}

\begin{document}

\title[Hybrid Normed Ideal Perturbations of $n$-tuples]{Hybrid Normed Ideal Perturbations of $n$-tuples of Operators I}
\author{Dan-Virgil Voiculescu}
\address{D.V. Voiculescu \\ Department of Mathematics \\ University of California at Berkeley \\ Berkeley, CA\ \ 94720-3840}
\thanks{Research supported in part by NSF Grant DMS-1665534.}
\date{}
\begin{abstract}
In hybrid normed ideal perturbations of $n$-tuples of operators, the normed ideal is allowed to vary with the component operators. We begin extending to this setting the machinery we developed for normed ideal perturbations based on the modulus of quasicentral approximation and an adaptation of our non-commutative generalization of the Weyl--von~Neumann theorem. For commuting $n$-tuples of hermitian operators, the modulus of quasicentral approximation remains essentially the same when $\cC_n^-$ is replaced by a hybrid $n$-tuple $\cC_{p_1,\dots}^-,\dots,\cC^-_{p_n}$, $p_1^{-1} + \dots + p_n^{-1} = 1$. The proof involves singular integrals of mixed homogeneity.
\end{abstract}

\maketitle

\begin{center}
Dedicated to Alain Connes on the occasion of his 70th birthday.
\end{center}

\bigskip
\section{Introduction}
\label{sec1}

In \cite{13} we adapted the Voiculescu non-commutative Weyl--von~Neumann type theorem \cite{12} to general normed ideals, which provided a new approach to normed ideal perturbation questions for $n$-tuples of operators on Hilbert space (see \cite{16} for a survey). A key quantity in this approach is the numerical invariant $k_{\Phi}(\tau)$, the modulus of quasicentral approximation where $\Phi$ is the norming function of a normed ideal $\cG_{\Phi}^{(0)}$ (see \cite{8} or \cite{11}) and $\tau = (T_1,\dots,T_n)$ is an $n$-tuple of operators. In particular, this technique has been successful in dealing with generalizations to commuting $n$-tuples of hermitian operators of the one hermitian operator Kuroda--Weyl--von~Neumann diagonalizability results and of the Kato--Rosenblum result on preservation of the Lebesgue absolutely continuous part under trace class perturbation (\cite{2}, \cite{7}, \cite{13}, \cite{14}). Since normed ideals provide the infinitesimals in non-commutative geometry there have also been technical uses of the machinery in non-commutative geometry (\cite{4}, \cite{5}, \cite{6}, \cite{9}).

The present paper, the first in a series, is the beginning of an extension of our approach to hybrid normed ideal perturbations. This means that instead of two $n$-tuples $\tau = (T_j)_{1 \le j \le n}$, $\tau' = (T'_j)_{1 \le j \le n}$ so that $T_j - T'_j \in \cG^{(0)}_{\Phi}$ for $1 \le j \le n$, we will consider now the more general situation $T_j - T'_j \in \cG_{\Phi_j}^{(0)}$ for $1 \le j \le n$, that is for each index $j$ there is now a different normed ideal $\cG_{\Phi_j}^{(0)}$ $(1 \le j \le n)$. Of course, much of the basics easily extends to the hybrid setting, with proofs that are almost verbatim repetitions of the same-ideal proofs. However, once we get to hybrid perturbations of $n$-tuples of commuting hermitian operators there is an unexpected new phenomenon: the threshold ideal $\cC_n^-$ (aka $\cC_{n,1}$ on the Lorentz scale) can be replaced by any hybrid $n$-tuple of normed ideals $(\cC_{p_1}^-,\dots,\cC_{p_n}^-)$ where $\sum_{1 \le j \le n} p_j^{-1} = 1$, $1 < p_j$ for $1 \le j \le n$. Indeed the generalized singular and absolutely continuous subspaces with respect to $(\cC_{p_1}^-,\dots,\cC_{p_n}^-)$ like in the case of $\cC_n^-$ coincide with the Lebesgue singular and Lebesgue absolutely continuous parts of the $n$-dimensional spectral measure. Actually even more, up to a factor of proportionality, the modulus of quasicentral approximation, is given by the same formula for all hybrid choices $(\cC_{p_1}^-,\dots,\cC_{p_n}^-)$ with $\sum_{1 \le j \le n} p_j^{-1} = 1$, which also includes the case when $p_1 = \dots = p_n = n$. While the general lines of the proofs are similar, one must refine some of the analysis involved with the singular integral underlying the non-vanishing of the modulus of quasicentral approximation, especially since one must replace the homogeneous singular integral with a new one with mixed homogeneity.

In addition to this introduction which is Section~\ref{sec1} and to the references, there are ten more sections.

Section~\ref{sec2} summarizes for the reader's convenience notation and definitions about normed ideals which we shall use.

Section~\ref{sec3} introduces the modulus of quasicentral approximation $k_{\varphi}(\tau)$ in the hybrid setting. Many of the simplest facts about $k_{\varphi}(\tau)$ are straightforward generalizations of \cite{13} and \cite{14} and the proofs could be omitted.

Section~\ref{sec4} deals with traces of sums of commutators, more precisely the commutators are between the operators in a given $n$-tuple and operators from the duals of the normed ideals. Traces of sums of such commutators are the source of lower bounds for the modulus of quasicentral approximation. This section is a straightforward generalization to the hybrid setting of results in \cite{15}.

Section~\ref{sec5} contains the hybrid version of the adaptation to normed ideals in \cite{13} of our main result in \cite{12} (see also \cite{1}). Again the proofs follow closely those in \cite{13} and were omitted.

Section~\ref{sec6} gives the hybrid generalization of the decomposition of the Hilbert space into $\Phi$-absolutely continuous and $\Phi$-singular subspaces with respect to an $n$-tuple of operators $\tau$ in \cite{14}. The straightforward generalizations of the proofs are omitted. In addition we provide also a strengthening of one of the results.

Section~\ref{sec7} deals with upper bounds for $k_{\varphi}$ and diagonalization $\mod \varphi$ results for $n$-tuples of commuting hermitian operators where $\varphi$ is a hybrid $n$-tuple of normed ideals. Upper bounds for $k_{\varphi}$ are much easier than lower bounds and this is why we can deal with larger classes of hybrid normed ideal perturbations and results are sharper. We obtain also diagonalization results in the case of spectral measures which are singular with respect to Lebesgue measure. These are, of course, the automatic consequence of vanishing of $k_{\varphi}$ results. The results we get in the case of $n = 2$, where the pair of norming functions are in duality, provide a large class of hybrid examples.

Sections~\ref{sec8} and \ref{sec9} deal with the estimate of the singular integral with mixed homogeneity which we need for the lower bound of $k_{\varphi}$ of an $n$-tuple of commuting hermitian operators. The general outline is similar to \cite{7}, but we need to refine some of the details and to use some stronger facts about absolutely convergent Fourier series.

In Section~\ref{sec10}, based on the previous technical facts, we prove the formula for $k_{\varphi}(\tau)$ where $\tau$ is an $n$-tuple of commuting hermitian operators and $\varphi$ a hybrid $n$-tuple of normed ideals $(\cC_{p_1}^-,\dots,\cC_{p_n}^-)$ with $p_1^{-1} + \dots + p_n^{-1} = 1$ and $1 \le p_j$, $1 \le j \le n$. In particular, this also implies that the $\varphi$-singular and $\varphi$-absolutely continuous subspaces are the same as the Lebesgue-singular and Lebesgue-absolutely continuous subspaces for the spectral measure of $\tau$.

In Section~\ref{sec11} we show that the commutants mod normed ideals have a natural generalization to the hybrid setting and we show that the relation between the modulus of quasicentral approximation and approximate units for the compact ideal \cite{17} carries over to this more general setting.

\section{Preliminaries}
\label{sec2}

Most of this section is about the notation we will use for operators and normed ideals of operators (\cite{8}, \cite{11}).

Let $\cH$ be a separable infinite-dimensional complex Hilbert space. By $\cL(\cH)$, $\cL_+(\cH)$, $\cK(\cH)$, $\cP(\cH)$, $\cR(\cH)$, $\cR_1^+(\cH)$ (or when $\cH$ is not in doubt simply $\cL$, $\cL_+$, $\cK$, $\cP$, $\cR$, $\cR_1^+$) we shall denote respectively the bounded operators on $\cH$, the positive bounded operators on $\cH$, the positive bounded operators on $\cH$, the compact operators, the finite rank hermitian projections, the finite rank operators and the finite rank positive contractions.

Let ${\hat c}$ be the space of sequence $(\xi_j)_{j \in \bbN}$, $\xi_j \in \bbR$ with finite support. A norming function $\Phi$ is a function on ${\hat c}$ taking values in $\bbR$ so that
\[
\begin{aligned}
\mbox{I.} &\quad\xi \ne 0 \Rightarrow \Phi(\xi) > 0 \\
\mbox{II.} &\quad\Phi(\alpha\xi) = |\alpha|\Phi(\xi),\ \alpha \in \bbR \\
\mbox{III.} &\quad\Phi(\xi+\eta) \le \Phi(\xi) + \Phi(\eta) \\
\mbox{IV.} &\quad\Phi((1,0,0,\dots)) = 1 \\
\mbox{V.} &\quad\Phi((\xi_j)_{j \in \bbN}) = \Phi((\xi_{\pi(j)})_{j \in \bbN}) \text{ if $\pi$ is a permutation of $\bbN$.}
\end{aligned}
\]
The set of norming functions will be denoted by $\cF$.

If $T \in \cR(\cH)$ and $\Phi \in \cF$ then
\[
|T|_{\Phi} = \Phi((s_j)_{j \in \bbN})
\]
where $(s_j)_{n \in \bbN}$ are the eigenvalues of $(T^*T)^{1/2}$ (multiple eigenvalues repeated according to multiplicity).

If $T \in \cL(\cH)$ it is a fact that we can extend the definition of $|\quad|_{\Phi}$ by putting
\[
|T|_{\Phi} = \sup_{P \in \cP} |TP|_{\Phi}
\]
which may now also take the value $+\infty$. The set $\cG_{\Phi}(\cH)$ or simply $\cG_{\Phi}$, when $\cH$ is not in doubt, is the set
\[
\{T \in \cL(\cH) \mid |T|_{\Phi} < \infty\}.
\]
Then $\cG_{\Phi}$ is an ideal in $\cL$ and a Banach space with respect to the norm $|\quad|_{\Phi}$. The closure of $\cR$ in $\cG_{\Phi}$ is denoted by $\cG_{\Phi}^{(0)}$ and is also an ideal in $\cL$. The norming function $\Phi$ is called mononorming if $\cG_{\Phi}^{(0)} = \cG_{\Phi}$ and binorming if $\cG_{\Phi}^{(0)} \ne \cG_{\Phi}$. If $T \in \cG_{\Phi}$ and $A,B \in \cL$ we have $|ATB|_{\Phi} \le \|A\||T|_{\Phi}\|B\|$. If rank $T = 1$ then $|T|_{\Phi} = \|T\|$. The norming functions
\[
\Phi_p((\xi_j)_{j \in \bbN}) = \left( \sum_{j \in \bbN} |\xi_j|^p\right)^{1/p}(1 \le p < \infty)
\]
and
\[
\Phi_{\infty}((\xi_j)_{j \in \bbN}) = \sup_{j \in \bbN}|\xi_j|
\]
define the Schatten--von~Neumann classes $\cC_p = \cG_{\Phi_p} = \cG_{\Phi_p}^{(0)}$ if $1 \le p < \infty$ and $\cG_{\Phi_{\infty}}^{(0)} = \cK$, $\cG_{\Phi_{\infty}} = \cL$, $|\quad|_{\infty} = \|\quad\|$ if $p = \infty$. We shall use the notation $|\quad|_p = |\quad|_{\Phi_p}$ if $1 \le p < \infty$ and we have $|T|_p = (\mbox{Tr}((T^*T)^{p/2})^{1/p}$.

If $\pi = (\pi_j)_{j \in \bbN}$ is a decreasing sequence $\pi_1 \ge \pi_2 \ge \dots$ of positive numbers so that $\lim_{j \to \infty} \pi_j = 0$ and $\sum_{j \in \bbN} \pi_j = \infty$ there are two norming functions $\Phi_{\pi}$ and $\Phi_{\pi^*}$ defined by
\[
\Phi_{\pi}((\xi_j)_{j \in \bbN}) = \sum_{j \in \bbN} \pi_j\xi^*,
\]
and
\[
\Phi_{\pi^*}((\xi_j)_{j \in \bbN}) = \sup_{n \in \bbN} \left( \left.\sum_{1 \le j \le n} \xi_j^*\right/ \sum_{1 \le j \le n} \pi_j\right)
\]
where $(\xi_j^*)_{j \in \bbN}$ is the decreasing rearrangement of $(|\xi_j|)_{j \in \bbN}$. Then $\Phi_{\pi}$ is mononorming, while $\Phi^*_{\pi}$ may not be. We shall denote $\cG_{\Phi_{\pi}} = \cG_{\Phi_{\pi}}^{(0)}$, $\cG_{\Phi_{\pi^*}}^{(0)}$, $\cG_{\Phi_{\pi^*}}$ by $\cG_{\pi}$, $\cG_{\pi^*}^{(0)}$, $\cG_{\pi^*}$, respectively. If $\pi_j = j^{-1+1/p}$ $(1 < p \le \infty)$, we shall denote $\cG_{\pi}$ by $\cC_p^-$ and its norm and norming function by $|\quad|_p^-$, $\Phi_p^-$ and if $\pi_j = j^{-1}/p$ $(1 \le p < \infty)$ we denote $\cG_{\pi^*}$ by $\cC_p^+$ and its norm and norming function by $|\quad|^+_p$, $\Phi_p^+$. We have $\cC_p^- \subset \cC_p \subset \cC_p^+$ and if $p < r$ then $\cC_p^+ \subset \cC_r^-$. The ideals $\cC^-_p$, $\cC_p^+$ actually correspond to $\cC_{p,1}$ and $\cC_{p,\infty}$ among the two-indices Lorentz-type $\cC_{p,q}$ ideals. The ideals $\cC_{\infty}^-$ and $\cC_1^+$ are often called Macaev ideals.

For every norming function $\Phi \in \cF$ there is a conjugate norming function $\Phi^* \in \cF$ so that $\cG_{\Phi^*}$ is the dual of $\cG_{\Phi}^{(0)}$ the duality pairing being
\[
\cG_{\Phi}^{(0)} \times \cG_{\Phi^*} \ni (X,Y) \to \mbox{Tr } XY \in \bbC.
\]
If $1/p + 1/q = 1$ then $\cC_q^+$ is the dual of $\cC_p^-$ $(1 < p \le \infty)$. Also in general $\cG_{\pi}$ is the dual of $\cG_{\pi^*}^{(0))}$ and $\cG_{\pi^*}$ is the dual of $\cG_{\pi}$.

We will use the following notation for $n$-tuples of operators. By $[n]$ we denote the set $\{1,\dots,n\}$ and by $\cL([n]\mid \cH)$ or just $\cL([n])$ we denote the maps $\tau: [n] \to \cL(\cH)$, which amounts to the same as $n$-tuples $(\tau(1),\dots,\tau(n))$ of operators. Similarly we use the notation $\cK([n] \mid \cH)$, $\cR([n] \mid \cH)$, $\cP([n] \mid \cH)$, $\cG_{\Phi}([n] \mid \cH)$, $\cG_{\Phi}^{(0)}([n] \mid \cH)$, etc., or just $\cK([n])$, $\cR([n])$, $\cP([n])$, $\cG_{\Phi}([n])$, $\cG_{\Phi}^{(0)}([n])$, etc., and $\|\tau\| = \max_{1 \le j \le n} \|\tau(j)\|$, $|\tau|_{\Phi} = \max_{1 \le j \le n} |\tau(j)|_{\Phi}$. Since we will deal with hybrid normed-ideal perturbations, we will also consider $\varphi \in \cF([n])$, $n$-tuples of norming functions, that is $\varphi: [n] \to \cF$, in which $\cG_{\varphi}^{(0)} \subset \cG_{\varphi}([n]) \subset \cL([n])$ are the subsets of $\cL([n])$ consisting of those $\tau$ so that $\tau(j) \in \cG_{\varphi(j)}^{(0)}$ and respectively $\tau(j) \in \cG_{\varphi(j)}$, $1 \le j \le n$ and $|\tau|_{\varphi} = \max_{1 \le j \le n} |\tau(j)|_{\varphi(j)}$. To distinguish between single norming functions and $n$-tuples of norming functions, we shall use capital greek letters like $\Phi$ for the singleton and lower case greek letters like $\varphi$ for the case of an $n$-tuple.

Further notation for operations on $n$-tuples are: if $\tau \in \cL([n])$ and $X,Y \in \cL$ then we define the $n$-tuples $\tau^*$, $X\tau Y$ and $[\tau,X]$ by
\[
\begin{aligned}
\tau^*(j) &= (\tau(j))^*, \\
X\tau Y(j) &= X\tau(j)Y \\
[\tau,X](j) &= [\tau(j),X].
\end{aligned}
\]
If $\sigma,\tau \in \cL([n])$ then $\sigma + \tau$ is defined by
\[
(\sigma+\tau)(j) = \sigma(j) + \tau(j).
\]
Also, if $\sigma \in \cL([m])$, $\tau \in \cL([n])$ then $(\sigma,\tau)$ is defined by
\[
\begin{aligned}
(\sigma,\tau)(i) &= \sigma(i) \mbox{ if $1 \le i \le m$} \\
(\sigma,\tau)(m+j) &= \tau(j) \mbox{ if $1 \le j \le n$.}
\end{aligned}
\]
Sometimes $(\sigma,\tau)$ will also be denoted $\sigma \amalg \tau$.

\section{The hybrid modulus of quasicentral approximation}
\label{sec3}

This section deals with the generalization to the hybrid setting of multiple norming functions of basic facts about the modulus of quasicentral approximation $k_{\Phi}(\tau)$. We shall often omit the obvious generalizations of the proofs from the case of one norming function and only indicate, when needed, where to find that proof.

If $\tau \in \cL([n] \mid \cH)$ and $\varphi \in \cF([n])$, we define
\[
k_{\varphi}(\tau) = \liminf_{A \in \cR_1^+(\cH)} |[\tau,A]|_{\varphi}
\]
where the liminf is with respect to the natural order on $\cR_1^+$. The key difference with the case of one norming function is that now
\[
|[\tau,A]|_{\varphi} = \max_{1 \le j \le n} |[\tau(j),A]|_{\varphi(j)}
\]
and the norming functions $\varphi(j)$, $1 \le j \le n$ depend on the value of $j$.

A $\tau$-quasicentral approximate unit relative to $\varphi$ is by definition a sequence $A_k \uparrow I$, $k \to \infty$, $A_k \in \cR_1^+$ so that
\[
\lim_{k \in \infty} |[\tau,A_k]|_{\varphi} = 0.
\]
The existence of a $\tau$-quasicentral approximate unit relative to $\varphi$ is equivalent to $k_{\varphi}(\tau) = 0$.

The proofs of the next three propositions are quite standard and easy.

\bigskip
\noindent
{\bf 3.1. Proposition.} {\em Let $\tau \in \cL([n])$, $\varphi \in \cF([n])$ and $\Phi \in \cF$ so that $\varphi(j) \le \Phi$ for all $1 \le j \le n$. Let further $A_i \in \cG_{\Phi}^{(0)}$, $0 \le A_i \le I$ where $i \in \cJ$, $\cJ$ being a directed set and assume $w - \lim_{i \in \cJ} A_i = I$. Then we have
\[
k_{\varphi}(\tau) \le \liminf_{i \in \cJ} |[\tau,A_i]|_{\varphi}.
\]
}

\bigskip
\noindent
{\bf 3.2. Proposition.} {\em Given $\tau \in \cL([n])$ and $\varphi \in \cF([n])$ there is a sequence $A_m \in \cR_1^+$, $m \in \bbN$, so that $A_m \uparrow I$ as $m \to \infty$ and
\[
k_{\varphi}(\tau) = \lim_{m \to \infty} |[\tau,A_m]|_{\varphi}.
\]
}

\bigskip
\noindent
{\bf 3.3. Proposition.} {\em Let $\tau \in \cL([n])$, $\varphi \in \cF([n])$ and $\sigma \in \cG_{\varphi}^{(0)}([n])$, then we have
\[
k_{\varphi}(\tau+\sigma) = k_{\varphi}(\tau).
\]
}

\bigskip
The next proposition is the generalization of Proposition~$1.5$ of \cite{15} and is proved along the same lines.

\bigskip
\noindent
{\bf 3.4. Proposition.} {\em Let $\tau \in \cL([n])$ and let $A_m \in \cK([n])$, $m \in \bbN$ be so that $A_m = A_m^*$ and $s - \lim_{m \to \infty} A_m = I$ and $[\tau,A_m] \in \cG_{\varphi}([n])$ for all $m \in \bbN$.  Then we have
\[
k_{\varphi}(\tau) \le \alpha \liminf_{m \to \infty} |[\tau,A_m]|_{\varphi}
\]
where $\alpha$ is a universal constant.}

\bigskip
The next proposition generalizes Proposition~1.4 of \cite{13} and its proof is omitted.

\bigskip
\noindent
{\bf 3.5. Proposition.} {\em Let $\tau^{(m)} \in \cL([n])$, $m \in \bbN$ and let $\varphi \in \cF([n])$. We have
\[
\begin{aligned}
\max_{1 \le j \le 2} k_{\varphi}(\tau^{(j)}) &\le k_{\varphi}(\tau^{(1)} \oplus \tau^{(2)}) \le k_{\varphi}(\tau^{(1)}) + k_{\varphi}(\tau^{(2)}) \\
k_{\varphi}\left(\bigoplus_{m \in \bbN} \tau^{(m)}\right) &= \lim_{m \to \infty} k_{\varphi}\left(\bigoplus_{j=1}^m \tau^{(j)}\right).
\end{aligned}
\]
}

\bigskip
A norming function $\Phi$ has property $(\Sigma)$ (\cite{14}) if
\[
\lim_{m \to \infty} \frac {1}{m} |X_1 \oplus \dots \oplus X_m|_{\Phi} = 0
\]
whenever $X_j \in \cR$, $\sup_{j \in \bbN} |X_j|_{\Phi} < \infty$. We shall say that $\varphi \in \cF([n])$ has property $(\Sigma)$ if the norming functions $\varphi(1),\dots,\varphi(n)$ have property $(\Sigma)$. The next proposition is an immediate generalization of Proposition~2.4 of \cite{14}, the proof of which being along the same lines will be omitted.

\bigskip
\noindent
{\bf 3.6. Proposition.} {\em Let $\tau \in \cL([n])$ and assume $\varphi \in \cF([n])$ has property $(\Sigma)$. Then $k_{\varphi}(\tau)$ is either $0$ or $\infty$.
}

\bigskip
The Schatten--von~Neumann classes $\cC_p$ with $1 < p$ have property $(\Sigma)$ and also $\cG_{\pi}$ when
\[
\lim_{m \to \infty} \left( \sup_{k \in \bbN} (\pi_k/m\pi_{km})\right) = 0
\]
has property $(\Sigma)$.

The next result is the hybrid version of Proposition~1.6 in \cite{13} and the proof being almost the same, will be omitted.

\bigskip
\noindent
{\bf 3.7. Proposition.} {\em Let $\tau^{(j)} \in \cL([n])$, $1 \le j \le m$ and let $c^{(j)} \in \cL([n])$ where $c^{(j)}(k) = c_k^{(j)}I$, $1 \le j \le m$, $1 \le k \le n$, $c_k^{(j)} \in \bbC$. Then we have
\[
k_{\varphi}(\tau^{(1)} \oplus \dots \oplus \tau^{(m)}) = k_{\varphi}((\tau^{(1)} - c^{(1)}) \oplus \dots \oplus (\tau^{(m)}-c^{(m)}))
\]
where $\varphi \in \cF([n])$.
}

\section{Traces of sums of commutators}
\label{sec4}

This section is a straightforward generalization of Section~II of \cite{15} and the proofs will be omitted. If $\Phi \in \cF$ is a norming function and $\Phi^*$ its conjugate let $p_{\Phi^*}: \cG_{\Phi^*} \to \cG_{\Phi^*}/\cG_{\Phi^*}^{(0)}$ be the canonical map and $|X|_{\Phi^*}^\sim$ will denote the norm of $p_{\Phi^*}(X)$ in $\cG_{\Phi^*}/\cG_{\Phi^*}^{(0)}$. If $\varphi \in \cF([n])$ and $\varphi^*$ is the conjugate $n$-tuple, we shall define $|\chi|^{\sim}_{\varphi^*} = \sum_{1 \le j \le n} |\chi(j)|^{\sim}_{\varphi^*}$, if $\chi \in \cG_{\varphi^*}([n])$. We have chosen here to combine the $n$ seminorms by summation instead of a max for duality reason. We also consider the canonical map 
\[
p_{\varphi^*}: \cG_{\varphi^*}([n]) \to \cG_{\varphi^*}([n])/\cG_{\varphi^*}^{(0)}([n])
\]
and also denote the last quotient by $(\cG/\cG^{(0)})_{\varphi^*}([n])$ which stands for $n$-tuples $(x(1),\dots,x(n))$ with $x(j) \in \cG_{\varphi^*(j)}/\cG_{\varphi^*(j)}^{(0)}$, $1 \le j \le n$.

\bigskip
\noindent
{\bf 4.1. Proposition.} {\em Let $\tau \in \cL([n])$, $\varphi \in \cF([n])$ and $\chi \in \cG_{\varphi^*}(n)$.
\begin{itemize}
\item[(i)] If $Y = \sum_{1 \le j \le n} [\tau(j),\chi(j)] \in \cC_1 + \cL_+$, then we have
\[
|\mbox{\rm Tr } Y| \le k_{\varphi}(\tau)|\chi|^{\sim}_{\varphi^*}.
\]
\item[(ii)] If $\tau$ is an $n$-tuple of unitary operators and
\[
Y = \sum_{1 \le j \le n} (\tau(j)\chi(j)\tau(j)^* - \chi(j)) \in \cC_1 + \cL_+
\]
then we have
\[
|\mbox{\rm Tr } Y| \le k_{\varphi}(\tau)|\chi|^{\sim}_{\varphi^*}\,.
\]
\end{itemize}
}

\bigskip
\noindent
{\bf 4.1. Corollary.} {\em Let $\tau \in \cL([n])$, $\varphi \in \cF([n])$, $\chi \in \cG_{\Phi^*}([n])$ and assume $k_{\varphi}(\tau) = 0$.
\begin{itemize}
\item[(i)] If $Y = \sum_{j \in [n]} [\tau(j),\chi(j)] \in \cC_1 + \cL_+$, then we have $Y \in \cC_1$ and $\mbox{\rm Tr } Y = 0$.
\item[(ii)] If the $\tau(j)$ are unitary, $1 \le j \le n$ and $Y = \sum_{j \in [n]} (\tau(j)\chi(n)\tau(j)^* - \chi(j)) \in \cC_1 + \cL_+$, then we have $Y \in \cC_1$ and $\mbox{\rm Tr } Y = 0$.
\end{itemize}
}

\bigskip
\noindent
{\bf 4.2. Proposition.} {\em Let $\tau \in \cL([n])$, $\varphi \in \cF([n])$ and assume $0 < C < k_{\varphi}(\tau)$.
\begin{itemize}
\item[(i)] If $\tau = \tau^*$, then there is $\chi \in \cG_{\varphi^*}([n])$, $\chi = \chi^*$, $\sum_{j \in [n]} |\chi(j)|_{\varphi^*(j)} \le 1$ so that
\[
Y = i \sum_{j \in [n]} [\chi(j),\tau(j)] \in \cC_1 + \cL_+
\]
and $\mbox{\rm Tr } Y \ge C$. Moreover, if $k_{\varphi}(\tau) < \infty$ then $Y \in \cC_1$.
\item[(ii)] If the $\tau(j)$ are unitary, $j \in [n]$ then there is $\chi \in \cG_{\varphi^*}([n])$, $\chi = \chi^*$, $\sum_{j \in [n]} |\chi(j)|_{\varphi^*(j)} \le 1$ so that
\[
Y = \sum_{j \in [n]} (\tau(j)\chi(j)\tau(j)^* - \chi(j)) \in \cC_1 + \cL_+
\]
and $\mbox{\rm Tr } Y \ge C$. Moreover, if $k_{\varphi}(\tau) < \infty$ then $Y \in \cC_1$.
\end{itemize}
}

\bigskip
To use the preceding result in case $\tau$ is not hermitian one need only to pass to $(2n)$-tuples and notice that
\[
\begin{aligned}
k_{\varphi}(\tau) &\ge k_{\varphi \amalg \varphi}\left( \frac {1}{2} (\tau+\tau^*) \amalg \frac {1}{2i} (\tau-\tau^*)\right) \\
&\ge \frac {1}{2} k_{\varphi}(\tau).
\end{aligned}
\]

\bigskip
\noindent
{\bf 4.3. Proposition.} {\em Let $\tau \in \cL([n])$, $\varphi \in \cF([n])$ be so that $k_{\varphi}(\tau) > 0$. Then there are decreasing sequences $\pi^{(k)} = (\pi_j^{(k)})_{j=1}^{\infty}$ with $\sum_{j \in \bbN} \pi_j^{(k)} = \infty$, $k \in [n]$, so that $\cG_{\psi} \supset \cG_{\varphi}^{(0)}$ and $k_{\psi}(\tau) > 0$ where $\psi \in \cF([n])$, $\psi(k) = \Phi_{\pi^{(k)}}$, $k \in [n]$.
}

\bigskip
The statement of this last proposition, the straightforward generalization of Proposition~2.6 in \cite{15}, includes also a small correction to that statement, the sequences $\pi^{(k)}$ may not converge to zero, in which case $\cG_{\pi(k)} = \cC_1$. We had overlooked the possibility of $\cC_1$.

\section{Adaptation of the Voiculescu non-commutative Weyl--von~Neumann type theorem for hybrid perturbations}
\label{sec5}

The canonical homomorphism to the Calkin algebra will be denoted by $p: \cL(\cH) \to \cL/\cK(\cH)$ or abbreviated $p: \cL \to \cL/\cK$. The next two theorems are extensions to the hybrid case of variations of Theorem~2.4 in \cite{13} and Corollary~2.5 in \cite{13}.

\bigskip
\noindent
{\bf 5.1. Theorem.} {\em Let $I \in \cA \subset \cL(\cH)$ be a $C^*$-subalgebra generated by an $n$-tuple $\tau([n])$ of operators and let $\rho: p(\cA) \to \cL(\cH)$ be a unital $*$-homomorphism so that $k_{\varphi}(\rho(p(\tau([n])))) = 0$ for a given $\varphi \in \cF([n])$. Then there are unitary operators $U_m: \cH \to \cH \oplus \cH$, $m \in \bbN$ so that
\[
U_m^*(\tau(k) \oplus \rho(p(\tau(k))))U_m - \tau(k) \in \cG_{\varphi(k)}^{(0)},\ k \in [n]
\]
and
\[
\lim_{m \to \infty} |U_m^*(\tau(k) \oplus \rho(p(\tau(k))))U_m - \tau(k)|_{\varphi(k)} = 0, k \in [n].
\]
}

\bigskip
\noindent
{\bf 5.2. Theorem.} {\em Let $\cA$ be a $C^*$-algebra with unit which is generated by a $n$-tuple $a([n])$ of elements and let $\varphi \in \cF([n])$. Let $\rho_j: \cA \to \cL(\cH)$, $j = 1,2$, be unital $*$-homomorphisms so that $\ker p\circ \rho_j = 0$, $j = 1,2$ and $k_{\varphi}(\rho_j(a([n]))) = 0$, $j = 1,2$. Then there are unitary operators $V_m: \cH \to \cH$, $m \in \bbN$, so that
\[
V_m\rho_1(a(k)) - \rho_2(a(k))V_m \in \cG_{\varphi(k)}^{(0)}, k \in [n]
\]
and
\[
\lim_{m \to \infty} |V_m\rho_1(a(k)) - \rho_2(a(k))V_m|_{\varphi(k)} = 0, k \in [n].
\]
}

\bigskip
The proofs are the obvious extensions to the hybrid case of the proofs in \cite{13} and will be omitted.

\section{The $\varphi$-singular and $\varphi$-absolutely continuous subspace of an $n$-tuple}
\label{sec6}

We begin with the extension to the hybrid case of the technical lemma underlying the definition of the $\Phi$-singular and $\Phi$-absolutely continuous subspaces for a $n$-tuple $\tau \in \cL([n])$, Lemma~$1.1$ in \cite{14}.

\bigskip
\noindent
{\bf 6.1. Lemma.} {\em Let $\tau \in \cL([n])$ and $\varphi \in \cF([n])$. We have $k_{\varphi}(\tau) = 0$ if and only if there is a sequence $A_m \in \cR$, $m \in \bbN$, $A_m \ge 0$ so that $A = w - \lim_{m \to \infty} A_m$ exists, $\lim_{m \to \infty} |[A_m,\tau]|_{\varphi} = 0$ and $\ker A = 0$.
}

\bigskip
\noindent
{\bf 6.1. Definition.} {\em Let $\tau \in \cL([n])$ and $\varphi \in \cF([n])$. We denote by $\cP_{\varphi}(\tau)$ the hermitian projections $P$ so that $[P,\tau] = 0$ and $k_{\varphi}(\tau \mid P\cH) = 0$. The $\varphi$-singular projection $E_{\varphi}^0(\tau)$ and $\varphi$-absolutely continuous projection $E_{\varphi}(\tau)$ of $\tau$ are defined to be
\[
\begin{aligned}
E_{\varphi}^0(\tau) &= \underset{P \in \cP_{\varphi}(\tau)}{\bigvee} P, \\
E_{\varphi}(\tau) &= I - E_{\varphi}^0(\tau).
\end{aligned}
\]
}

\bigskip
The basic fact about such subspaces Theorem~$1.2$ of \cite{14} also immediately extends.

\bigskip
\noindent
{\bf 6.1. Theorem.} {\em Let $\tau \in \cL([n])$ and $\varphi \in \cF([n])$.}

a) {\em We have $E_{\varphi}^0(\tau) \in \cP_{\varphi}(\tau)$, in particular $k_{\varphi}(\tau \mid E_{\varphi}^0(\tau)\cH) = 0$. Moreover $E_{\varphi}^0(\tau)$ is in $Z(W^*(\tau([n])))$, the center of the von~Neumann algebra of $\tau([n])$ and hence also $E_{\varphi}(\tau) \in Z(W^*(\tau([n])))$.}

b) {\em If $\Phi \in \cF$ is so that $\Phi \ge \varphi(k)$, $k \in [n]$ and $A_m = A_m^* \in \cG_{\Phi}^{(0)}$, $m \in \bbN$ are so that
\[
\sup_{m \in \bbN} \|A_m\| < \infty \mbox{ and } \lim_{m \to \infty} |[\tau,A_m]|_{\varphi} = 0
\]
then we have $s = \lim_{m \to \infty} A_mE_{\varphi}(\tau) = 0$.
}

\bigskip
Actually part b) of the preceding theorem can easily be given a stronger and more convenient form.

\bigskip
\noindent
{\bf 6.1. Proposition.} {\em Let $\tau \in \cL([n])$ and $\varphi \in \cF([n])$. If $A_m = A_m^* \in \cK$, $m \in \bbN$ are so that
\[
\sup_{m \in \bbN} \|A_m\| < \infty \mbox{ and } \lim_{m \to \infty} |[\tau,A_m]|_{\varphi} = 0
\]
then we have $s - \lim_{n \to \infty} A_mE_{\varphi}(\tau) = 0$.
}

\bigskip
\noindent
{\bf {\em Proof.}} Let $C > \|A_m\|$ for all $m \in \bbN$ and $g: (-C-1,C+1) \to \bbR$ a $C^{\infty}$-function such that
\[
\mbox{supp } g \subset [-C-1/2,-\epsilon/2] \cup [\epsilon/2,C+1/2],
\]
$|g(t)| \le |t|$ for all $t \in (-C-1,C+1)$, and $g(t) = t$ if $t \in [-C,-\epsilon] \cup [\epsilon,C]$. Then $g(A_m) \in \cR$, $\sup_{m \in \bbN} \|g(A_m)\| < \infty$, $\|g(A_m) - A_m\| \le \epsilon$, and $\lim_{m \to \infty} |[g(A_m),\tau]|_{\varphi} = 0$. By Theorem~$6.1$~b) we infer that
\[
s - \lim_{m \to \infty} g(A_m)E_{\varphi}(\tau) = 0.
\]
This implies that if $\xi \in E_{\varphi}(\tau)\cH$ we have
\[
\limsup_{m \to \infty} \|A_m\xi\| \le \epsilon\|\xi\| + \limsup_{m \to \infty} \|g(A_m)\xi\| = \epsilon\|\xi\|.
\]
Since this holds for all $\epsilon > 0$ we get the desired conclusion.\hfill$\qed$

\section{Upper bounds on $k_{\varphi}$ and diagonalization $\mod \varphi$ for $n$-tuples of commuting hermitian operators}
\label{sec7}

We begin with the consequence of our adapted Weyl--von~Neumann type Theorem~5.2 which relates diagonalization $\mod \varphi$ of a commuting $n$-tuple of hermitian operators and $k_{\varphi}$ (the hybrid extension of Corollary~2.6 in \cite{13}). The proof being along the same lines will be omitted.

\bigskip
\noindent
{\bf 7.1. Proposition.} {\em Let $\tau \in \cL([n])$ be an $n$-tuple of commuting hermitian operators and let $\varphi \in \cF([n])$. Then the following are equivalent:}

(i) $k_{\varphi}(\tau) = 0$.

(ii) {\em there is an $n$-tuple of hermitian operators $\delta \in \cL([n])$ which are simultaneously diagonal in some orthonormal basis, so that $\delta(k) - \tau(k) \in \cG_{\varphi(k)}^{(0)}$, $k \in [n]$ and $|\delta(k)-\tau(k)|_{\varphi(k)} < \epsilon$ for some given $\epsilon > 0$.
}

\bigskip
We would like to remark that the same result holds for unitary or normal operators. Since one can easily pass from one to the other, this being the case of commutative $C^*$-algebras, we shall deal with hermitian operators in what follows.

\bigskip
\noindent
{\b 7.2. Proposition.} {\em Let $\tau \in \cL([n])$ be an $n$-tuple of commuting hermitian operators and let $\varphi \in \cF([n])$ be so that
\[
\lim_{m \to \infty} m^{-1} \prod^n_{k=1} \varphi(k) (\underset{m}{\underbrace{1,\dots,1}},0,\dots) = B \mbox{ where } 0 < B < \infty.
\]
Then if $\tau$ has a cyclic vector, there is a constant $C$, which depends only on $\varphi$, so that
\[
k_{\varphi}(\tau) \le C\|\tau\|.
\]
}

\bigskip
\noindent
{\bf {\em Proof.}} It suffices to prove that $k_{\varphi}(\tau) \le C$, where $C$ depends only on $\varphi$ when $\|\tau\| < 1/2$. Replacing $(T_1,\dots,T_n)$ by $(T_1+1/2I,\dots,T_n+1/2I)$ which does not change $k_{\varphi}(\tau)$ we may assume $\sigma(\tau) \subset (0,1)^n$ instead of $\|\tau\| \le 1/2$. Let further $\xi$ denote the cyclic vector of $\tau$ and let
\[
N_k(m) = [\varphi(k)(\underset{m}{\underbrace{1,\dots,1}},0,\dots)]
\]
where $[\cdot]$ denotes the integer part. Then $\lim_{m \to \infty} m^{-1} N_1(m) \dots N_m(m) = B$. If $k_j \in \{1,\dots,N_j(m)\}$, $1 \le j \le n$, let $E(k_1,\dots,k_n)$ denote the spectral projection
\[
E\left( \tau; \prod_{j=1}^n [(k_j-1)N_j(m)^{-1},k_jN_j(m)^{-1}) \right)
\]
and $\xi(k_1,\dots,k_n) = E(k_1,\dots,k_n)\xi$. Then if $P_m$ is the orthogonal projection onto the subspace
\[
\cH_m = \underset{\begin{matrix}
{}^{1 \le k_j \le N_j(m)} \\
{}^{1 \le j \le n}
\end{matrix}}{\bigoplus} \bbC\xi(k_1,\dots,k_n)
\]
we have $[P_m,E(k_1,\dots,k_n)] = 0$ and $\underset{\begin{matrix}
{}^{1 \le k_j \le N_j(m)} \\
{}^{1 \le j \le n}
\end{matrix}}{\sum} E(k_1,\dots,k_m) = I$. Since both $T_j$ and $P_m$ are orthogonal direct sums of operators on the $E(k_1,\dots,k_n)\cH$, so is their commutator and since
\[
E(k_1,\dots,k_n)P_m = P_{\bbC\xi(k_1,\dots,k_n)}
\]
these are rank $\le 2$ and their norm is $\le 2/N_j(m)$. Thus we get
\[
\begin{aligned}
|[T_j,P_m]|_{\varphi(j)} &\le \frac {4}{N_j(m)} \varphi(j)(\underset{N_1(m)\dots N_n(m)}{\underbrace{1,\dots,1}},0,\dots) \\
&\le 4N_j(m)^{-1} N_j(N_1(m)\dots N_n(m)).
\end{aligned}
\]
Since $N_1(m)\dots N_n(m) \le (B+1)m$ for $m$ large enough, this gives $\limsup_{m \to \infty} |[T_j,P_m]_{\varphi(j)} \le 4(B+1)$, which proves our assertion.

\bigskip
\noindent
{\bf 7.3. Proposition.} {\em Let $\varphi \in \cF([n])$ be so that for all $n$-tuples $\tau$ of commuting hermitian operators which have a cyclic vector we have $k_{\varphi}(\tau) < C\|\tau\|$. Then, if $\tau$ is an $n$-tuple of commuting hermitian operators, the spectral measure of which is singular with respect to $n$-dimensional Lebesgue measure, then $k_{\varphi}(\tau) = 0$.
}

\bigskip
\noindent
{\bf {\em Proof.}} Using Proposition~3.5 it is easy to see that it suffices to prove the assertion under the additional assumption that $\tau$ has a cyclic vector. Thus we may assume $\tau$ is the $n$-tuple of multiplication operators by the coordinate functions in $L^2(\bbR^n,d\nu)$ where $\nu$ is a probability measure with compact support which is singular with respect to Lebesgue measure. Then we can find sets $\omega_m$ which are finite disjoint unions of sets
\[
[j_12^{-k},(j_1+1)2^{-k}) \times \dots \times [j_n2^{-k},(j_n+1)2^{-k})
\]
where $k$ depends only on $m$, and so that $\nu(\omega_m) \to 1$, while $\lambda_n(\omega_m) \to 0$ as $m \to \infty$. Since $\lambda_n(\omega_m) \to 0$ we can shift the cubes in $\omega_m$ to other disjoint positions, so that the new sets $\omega'_m$ have diameters shrinking to zero. If $\tau^{(m)}$ and $\tau'{}^{(m)}$ denote the $n$-tuples of multiplication operators in $L^2(\omega_m,\nu\mid \omega_m)$ and, respectively, in $L^2(\omega'_m,p_{m^*}|\nu|\omega_m|)$ where $p_m: \omega_m \to \omega'_m$ is the bijection resulting from shifting the cubes. The shifting of cubes and push-forward of measure which relates $\tau^{(m)}$ and $\tau'{}^{(m)}$ is precisely the situation to which we can apply Proposition~3.7, which gives $k_{\varphi}(\tau^{(m)}) = k_{\varphi}(\tau'{}^{(m)})$. On the other hand that $\lim_{m \to \infty} k_{\varphi}(\tau^{(m)}) = k_{\varphi}(\tau)$ is a consequence of the second assertion of Proposition~3.5. Actually since the diameter of $\omega'_m$ converges to zero an additional shift will insure that $\|\tau'{}^{(m)}\| \to 0$ as $m \to \infty$. Then observing that $\tau'{}^{(m)}$ also has a cyclic vector $k_{\varphi}(\tau^{(m)}) = k_{\varphi}(\tau'{}^{(m)}) \le C\|\tau'{}^{(m)}\|$ will give $k_{\varphi}(\tau^{(m)}) \to 0$ and hence $k_{\varphi}(\tau) = 0$.\hfill\qed

\bigskip
If $\Phi^-_p$ is the norming function for $\cC_p^-$ then $\Phi_p^-(\underset{m}{\underbrace{1,\dots,1}},1,0,\dots) = 1^{-1+1/p} + \dots + m^{-1+1/p}$ so that $\lim_{m \to \infty} m^{-1/p}\Phi^-_p(\underset{m}{\underbrace{1,\dots,1}},0,\dots) = 1/p$. Then if $p_j > 1$, $1 \le j \le n$ are so that $\sum_{1 \le j \le n} p_j^{-1} = 1$ we see that $\varphi \in \cF([n])$ defined by $\varphi(j) = \Phi_{p_j}$ satisfies the assumption of Proposition~7.2. Then in turn the conclusion of Proposition~7.2 means that the assumption of Proposition~7.3 is satisfied. Putting all this together we see that the following proposition holds.

\bigskip
\noindent
{\bf 7.4. Proposition.} {\em Let $p_j > 1$, $1 \le j \le n$ be so that $\sum_{1 \le j \le n} p_j^{-1}$ and let $\varphi \in \cF([n])$ be defined by $\varphi(j) = \Phi_{p_j}^-$, $1 \le j \le n$. Let further $\tau \in \cL([n])$ be an $n$-tuple of commuting hermitian operators. Then, if the spectral measure of $\tau$ is singular with respect to Lebesgue measure, then $k_{\varphi}(\tau) = 0$. Also, there is a constant $C$ which depends on $p_1,\dots,p_n$ so that $k_{\varphi}(\tau) \le C\|\tau\|$ if $\tau$ has a cyclic vector.
}

\bigskip
If $n = 2$ we can use Proposition~7.2 there is a quite general class of $\varphi \in \cF([z])$ to which it applies.

\bigskip
\noindent
{\bf 7.5. Proposition.} {\em Let $\pi = (\pi_j)_{j \in \bbN}$ be a decreasing sequence of positive numbers converging to zero and so that $\sum_{j \in \bbN} \pi_j = \infty$. Let $\varphi \in \cF([2])$ be defined to be $\varphi(1) = \Phi_{\pi}$, $\varphi(2) = \Phi^*_{\pi}$. Let further $\tau \in \cL([2])$ be a pair of commuting hermitian operators. Then there is a constant $C$ depending only on $\pi$, so that if $\tau$ has a cyclic vector, then $k_{\varphi}(\tau) \le C\|\tau\|$. Also, if the spectral measure of $\tau$ is singular with respect to Lebesgue measure then $k_{\varphi}(\tau) = 0$.
}

\bigskip
\noindent
{\bf {\em Proof.}} In order to apply Proposition~7.2 and Proposition~7.3 we will need only to check that the assumption of Proposition~7.2 is satisfied by $\varphi$. Indeed, we have 
\[
\Phi_{\pi}(\underset{m}{\underbrace{1,\dots,1}},0,\dots) = \pi_1 + \dots + \pi_m
\]
and 
\[
\Phi_{\pi^*}(\underset{m}{\underbrace{1,\dots,1}},0,\dots) = \max_{1 \le k \le m} \frac {k}{\pi_1+\dots+\pi_2}.
\]
Since $\pi$ is decreasing, the max in the right hand side of the last equality is actually $\frac {m}{\pi_1+\dots+\pi_m}$ and it follows that
\[
\Phi_{\pi}(\underset{m}{\underbrace{1,\dots,1}},0,\dots) \cdot \Phi_{\pi^*}(\underset{m}{\underbrace{1,\dots,1}},0,\dots) = m.
\]
\hfill\qed

\bigskip
\noindent
{\bf 7.6. Proposition.} {\em Let $p_j > 1$, $1 \le j \le n$ be so that $\sum_{1 \le j \le n} p_j^{-1}$ and $\varphi \in \cF([n])$ be so that $\varphi(j) = \Phi_{p_j}$, $1 \le j \le n$. Then if $\tau \in \cL([n])$ is an $n$-tuple of commuting hermitian operators then $k_{\varphi}(\tau) = 0$.
}

\bigskip
\noindent
{\bf {\em Proof.}} It suffices to prove the conclusion when $\tau$ has a cyclic vector. Then by Proposition~7.4 since $\Phi^-_p \ge \Phi_p$ we must have $k_{\varphi}(\tau) < \infty$. But $\varphi$ has property $(\Sigma)$ so by Proposition~3.6 we must have $k_{\varphi}(\tau) = 0$.\hfill\qed

\section{Absolutely convergent Fourier series in preparation for the lower bound}
\label{sec8}

Let $p_j > 1$, $1 \le j \le m$ be so that $\sum_{1 \le j \le m} p_j^{-1} = 1$ and let $\varphi \in \cF([m])$ be given by $\varphi(j) = \Phi_{p_j}^-$. In this section and in the next section we estimate certain singular integrals which using traces of commutators will show that $k_{\varphi}(\tau) > 0$ when $\tau$ is a commuting $m$-tuple of hermitian operators with Lebesgue absolutely continuous spectral measure. This will be achieved along broad lines similar to \cite{7} but the details will require several refinements for the hybrid case we consider. The kernels for the operators in $L^2$ of a cube with respect to $m$-dimensional Lebesgue measure will be of the form
\[
\mbox{sign}(x_j-y_j)|x_j-y_j|^{p_j-1}\left( \sum_{1 \le k \le m} |x_k-y_k|^{p_k}\right)^{-1}.
\]
We will fix notations choosing $j=1$ and instead of $\mbox{sign}(x_j-y_j)|x_j-y_j|^{p_j-1}$ we will deal with $(x_j-y_j)^{p_j-1}_+$ where $(x)_+ = \max(x,0)$. Our goal will be to show at the end of Section~\ref{sec9} that this operator is in $\cC_{q_j}^+$ where $p_j^{-1} + q_j^{-1} = 1$. In the present section we have collected some absolutely convergent Fourier series lemmata needed in the next section.

\bigskip
\noindent
{\bf 8.1. Lemma.} {\em Let $L > 0$ and let $\varphi$ be a $C^{\infty}$-function on $\bbR$ with $\mbox{supp } \varphi \subset (-L,L)$. Then, on $[-L,L]$ the function $f(x) = (x)_+^{\epsilon} \varphi(x)$, where $\epsilon > 0$, has absolutely convergent Fourier series. 
}

\bigskip
\noindent
{\bf {\em Proof.}} Since $(x)_+^{\epsilon}$ is increasing on $[-L,L]$ it has bounded variation and then also $(x)^{\epsilon}_+ \varphi(x)$ will have bounded variation. Moreover $f(x)$ is $\alpha$-Lipschitz for some $\alpha > 0$, since its derivative is in $L^p$ if $p(1-\epsilon) < 1$, $p > 1$. The absolute convergence of the Fourier series of $f$ follows then from Theorem~3.6 in Chapter~VI of \cite{18}.\hfill\qed

\bigskip
\noindent
{\bf 8.2. Lemma.} {\em Let $\psi: \bbR^m \to \bbR$ with $\mbox{supp } \psi \subset (-L,L)^m$ and which has absolutely convergent Fourier series in $[-L,L]^m$. Then $g(x_1,\dots,x_m) = (x_1)_+^{\epsilon}\psi(x_1,\dots,x_m)$, where $\epsilon > 0$, has absolutely convergent Fourier series in $[-L,L]^m$.
}

\bigskip
\noindent
{\bf {\em Proof.}} This is immediate from the preceding lemma. Let $0 < L_1 < L$ be so that $\mbox{supp } \psi \subset [-L_1,L_1]^m$ and $\varphi|[-L_1,L_1] \equiv 1$. Then $(x_1)_+^{\epsilon}\varphi(x_1)$ as a function on $[-L,L]^m$ has absolutely convergent Fourier series and \linebreak $(x_1)_+^{\epsilon}\psi(x_1,\dots,x_m) = ((x_1)_+^{\epsilon}\varphi(x_1))\psi(x_1,\dots,x_m)$.\hfill\qed

\bigskip
\noindent
{\bf 8.3. Lemma.} {\em Assume $p_j > 1$, $1 \le j \le m$ and $0 < r < R$ and $\psi$ be a $C^{\infty}$-function on $\bbR^m$ so that $\mbox{supp } \psi \subset ((-R,R)^m\backslash [-r,r]^m)$. Then $(|x_1|^{p_1} + \dots + |x_m|^{p_m})^{-1}\psi(x_1,\dots,x_m)$ has absolutely convergent Fourier series on $[-R,R]^m$. 
}

\bigskip
\noindent
{\bf {\em Proof.}} Let $\psi_1$ be another $C^{\infty}$-function with $\mbox{supp } \psi_1 \subset ((-R,R)^m\backslash [-r,r]^m)$ and so that $\psi_1 \mid \mbox{supp } \psi \equiv 1$. Then
\[
g(x_1,\dots,x_m) = (|x_1|^{p_1} + \dots + |x_m|^{p_m})\psi_1(x_1,\dots,x_m)
\]
has absolutely convergent Fourier series in $[-R,R]^m$. This follows by the same argument as in Lemma~8.1, though it follows also from weaker results on absolute convergence of Fourier series since $p_j > 1$ implies that $|x_j|^{p_j}$ is $C^1$. The fact that $g^{-1}\psi$ has absolutely convergent Fourier, after noticing that $g$ is bounded away from $0$ on $\mbox{supp } \psi$ is then a consequence of the analogue of Wiener's Lemma for $A(\bbT^m)$. (The case $m = 1$ of this result is Lemma~6.1 in Chapter~VIII of \cite{10}. The case of general $m$ can be found in \cite{3} Chapter~I, \S 5, Corollary~2 of Proposition~4 after using Chapter~I \S 6, Proposition~12 and Chapter~II, \S 1, Proposition~5 to see that the Banach algebra $L^1(\bbZ^m)$ is regular and has no radical.)

\section{$\cC_q^+$ estimates of the singular integrals for the lower bound}
\label{sec9}

In this section, in continuation of the preceding section, we prove the estimates on the singular integrals which will be used for the lower bound on $k_{\varphi}$. Overall, the method is similar to that which was used in \cite{7}, but sometimes the details need to be dealt differently, because the singular integrals which we consider have mixed homogeneity.

Let $\pi = (\pi_j)_{j \in \bbN}$ denote a sequence which decreases to zero and let $(\cG_{\pi^*},|\ |_{\pi^*})$ be the largest normed ideal corresponding to the norming function $\Phi_{\pi^*}(\xi) = \sup_{n \ge 1} (\sum_{1 \le k \le n} \xi_k^*/\sum_{1 \le k \le n} \pi_k)$. The sequence $\pi$ is regular if $\pi_1 + \dots + \pi_n \le C_0 n\pi_n$, in which case $T \in \cG_{\pi^*}$ iff $s_n = O(\pi_n)$, where the $s$-numbers $s_n$ are the eigenvalues of $(T^*T)^{1/2}$. In particular, if $\pi_j = j^{-1+1/p}$ where $1 < p < \infty$, we have $\cG_{\pi^*} = \cC_q^+$, where $q$ is the conjugate exponent to $p$, that is $p^{-1} + q^{-1} = 1$. The next proposition is Lemma~2.1 of \cite{7}, where the parameter $p$ has been specialized to $p = 1$ and the reader is referred to \cite{7} for the proof.

\bigskip
\noindent
{\bf 9.1. Proposition.} {\em Let $\pi$ be a regular sequence and let $T$ be a compact operator. Suppose $T$ has a decomposition
\[
T = \sum_{n \ge 1} \sum_{j \in I(n)} T_{n,j} \eqno (9.1)
\]
where
\[
T_{n,i}T_{n,j} = T^*_{n,i}T_{n,j} = T_{n,i}T_{n,j}^* = 0 \eqno (9.2)
\]
whenever $i \ne j$. Also, suppose that there are sequences $\beta_n$ and $\gamma_n$ and a constant $C > 0$ such that
\[
\gamma_n \uparrow +\infty \mbox{ with } \gamma_{n+1} \le C\gamma_n \eqno (9.3)
\]
\[
\beta_n \downarrow 0 \mbox{ with } \sum_{k > n} \beta_k \le C\beta_n \mbox{ for all $n$} \eqno (9.4)
\]
but
\[
\sum_{k \le n} \beta_k\gamma_k \le C\beta_n\gamma_n,\ \beta_n \le C\pi_{\gamma_n} \eqno (9.5)
\]
\[
|I(n)| \le \gamma_n \mbox{ for all } n \eqno (9.6)
\]
and
\[
|T_{n,j}| \le \beta_n \mbox{ for all $n$ and } j \in I(n). \eqno (9.7)
\]
Then $T \in \cG_{\pi^+}$ and $|T|_{\pi^*}$ depends only on $C$ and $C_0$.
}

\bigskip
The proposition will be used to prove that certain kernels give rise to operators in certain $\cC^+_q$ in $L^2(\mu)$, where $\mu$ is Lebesgue measure on $[-1,1]^m \subset \bbR^m$ and the kernel is defined using $K: \bbR^m\backslash \{0\} \to \bbR$ where
\[
K(x) = \frac {|x_1|^{p_1-1}\mbox{ sign } x_1}{|x_1|^{p_1} + \dots + |x_1|^{p_m}}\,,
\]
and $p_j > 1$, $1 \le j \le m$ are so that $\sum_{1 \le j \le m} p_j^{-1} = 1$. It will be more convenient to consider
\[
K_0(x) = (x_1)^{p_1-1}_+ (|x_1|^{p_1} + \dots + |x_m|^{p_m})^{-1}.
\]
The operator we shall study is $Tf(x) = \int K(x-y)d\mu(y)$ for $f \in L^2(\bbR^m,d\mu)$ or $T_0$ defined in the same way with $K_0$ instead of $K$. It is easily seen that $T$ is a linear combination of operators that are unitarily equivalent to $T_0$. Thus it will suffice to focus on $T_0,K_0$. {\em Because of this we shall change notation and denote by $T$ and $K$ the above $T_0$ and $K_0$} (and forget the initial $T,K$ for a while). Note also that $d\mu = \chi_{[-1,1]^m}d\lambda$ ($\chi_{\omega}$ denotes the indicator function of the set $\omega$).

Roughly our next step is to adapt the argument in the proof of Theorem~3.1 of \cite{7} to deal with the kernel with mixed homogeneity defined by
\[
K(x) = (x_1)_+^{p_1-1}(|x_1|^{p_1} + \dots + |x_m|^{p_m})^{-1}.
\]

It will be convenient to consider on $\bbR^m$ the function $\delta: \bbR^m \to [0,\infty)$ where $\delta(x_1,\dots,x_m) = \max(|x_1|^{p_1},\dots,|x_m|^{p_m})$ and to define
\[
D(x,r) = \{y \in \bbR^m \mid \delta(y-x) \le r\}
\]
when $x \in \bbR^m$ and $r \ge 0$. We shall use the dilations $\sigma(t): \bbR^m \to \bbR^m$, $t > 0$ given by
\[
\sigma(t)(x_1,\dots,x_m) = (t^{1/p_1}x_1,\dots,t^{1/p_m}x_m)
\]
so that $\sigma(t')\sigma(t'') = \sigma(t't'')$ and $\delta(\sigma(t)x) = t\delta(x)$, $\sigma(t)D(0,r) = D(0,tr)$, in particular $\sigma(r)D(0,1) = D(0,r)$ where actually $D(0,1) = [-1,1]^m$. We also have
\[
K(\sigma(t)x) = t^{-1/p_1}K(x).
\]
Remark also that
\[
\mu(D(x,r)) \le \lambda(D(x,r)) = \lambda(D(0,r)) = \lambda(\sigma(r)D(0,1)) = 2^mr
\]
where we used that $p_1^{-1} + \dots + p_m^{-1} = 1$.

By $B(x,r)$ we shall denote the closed ball of radius $r$ centered at $x \in \bbR^m$ with respect to the $l^{\infty}$-norm $\|x\| = \max_{1 \le j \le m} |x_j|$. Then, we have $D(0,1) = B(0,1) = [-1,1]^m$ and if $r \ge 1$ we have $D(x,r) \subset B(x,r)$.

\bigskip
\noindent
{\bf 9.2. Proposition.} {\em Let $K: \bbR^m\backslash\{0\} \to [0,\infty)$ be given by $K(x) = (x_1)^{p_1-1}_+ (|x_1|^{p_1} + \dots + |x_m|^{p_m})^{-1}$ and let $(Tf)(x) = \int K(x-y)d\mu(y)$ be the operator in $L^2(\mu)$ where $\mu$ is Lebesgue measure restricted to $[-1,1]^m = D(0,1)$. Then $T \in \cC_{q_1}^+$ where $q_1 = p_1(p_1-1)^{-1}$.
}

\bigskip
\noindent
{\bf {\em Proof.}} Since $\cC_{q_1}^+$ is $\cG_{\pi^*}$ for the choice of $\pi_j = j^{-1+1/p_1}$, we shall denote by $\pi_j$ these numbers throughout the proof and use Proposition~9.1 to show that $T \in \cG_{\pi^*}$.

\bigskip
We pick a $C^{\infty}$-function $\varphi_0$ with support in $D(0,2)$ and so that $\varphi_0 \mid D(0,1/2) \equiv 1$. We then define $\varphi(u) = \varphi_0(u) - \varphi_0(\sigma(2)u)$ and then $K_n(u) = K(u)\varphi (\sigma(2^n)u)$ and
\[
(S_nf)(x) = \int K_n(x-y)f(y)d\mu(y)
\]
which is an operator in $L^2(\mu)$. In a formal sense $T = \sum_{n \in \bbZ} S_n$.

Remark that $\varphi \mid D(0,1/4) = 0$ and there is $n_0 < 0$ so that $2D(0,1) \subset D(0,2^{-n_0})$ (the number $n_0$ depends on $p_1,\dots,p_m$). Then, if $(x,y) \in (\mbox{supp } \mu)^2 = (D(0,1))^2$ we have $\sigma(2^n)(x-y) \in D(0,1/4)$ if $n < n_0-2$ and hence $\varphi(\sigma(2^n)(x-y)) = 0$, $K_n(x-y) = 0$ and $S_n = 0$. Also, since $\sum_n \varphi(\sigma(2^n)u) = 1$ we get $\sum_n K_n = K$ which gives (in a formal sense) $T = \sum_n S_n = \sum_{n \ge n_0-2} S_n$.

We shall decompose each $S_n$ into parts that are ``localized''. Let $\theta \in C^{\infty}(\bbR^n)$ be a function with support in $B(0,10^{-1})$ and such that $\sum_{l \in \bbZ^m} \theta(x - 10^{-2}l) \equiv 1$.

For each $n \ge n_0 - 2$, we write $S_n = \sum_{l \in \bbZ^m} S_{n,l}$. Here 
\[
\begin{aligned}
S_{n,l}f(x) &= \theta(\sigma(2^n)x - 10^{-2}l)S_nf(x) \\
&= \int K_{n,l}(x,y)f(y)d\mu(y)
\end{aligned}
\]
where $K_{n,l}(x,y) = K(x-y)\varphi(\sigma(2^n)(x-y))\theta(\sigma(2^n)x-10^{-2}l)$. Since $K(x) = 2^{n/p_1}K(\sigma(2^n)x)$, we have
\[
K_n(u) = 2^{n/p_1}K(\sigma(2^n)u)\varphi(\sigma(2^n)u)
\]
and
\[
K_{n,l}(x,y) = 2^{n/p_1}K(\sigma(2^n)(x-y))\varphi(\sigma(2^n)(x-y))\theta(\sigma(2^n)x - 10^{-2}l).
\]
This implies that
\[
K_{n,l}(x,y) = 2^{n/p_1}K_{0,l}(\sigma(2^n)x,\sigma(2^n)y).
\]
Our next task will be to estimate the trace-norm (i.e., the $1$-norm) of $S_{n,l}$. We begin with a lemma.

\bigskip
\noindent
{\bf 9.1. Lemma.} {\em The operator defined in $L^2(\bbR^m,d\lambda)$ by $K_{0,0}(x,y)$ is in $\cC_1$.
}

\bigskip
\noindent
{\bf {\em Proof.}} Choose $L > 0$ so that $\mbox{supp } \varphi K \subset \left( -\frac {L}{4},\frac {L}{4}\right)^m$ and $\mbox{supp } \theta \subset \left( -\frac {L}{4},\frac {L}{4}\right)^m$. Let $e_{\alpha}(x) = \exp(\pi i \alpha \cdot xL^{-1})$ where $x \in [-L,L]^m$ and $\alpha \in \bbZ^m$. Then by Lemma~8.2 and Lemma~8.3 we have that
\[
\varphi(x)K(x) = \sum_{\alpha \in \bbZ^m} c_{\alpha}e_{\alpha}(x) \mbox{ with } \sum_{\alpha \in \bbZ^m} |c_{\alpha}| < \infty
\]
and clearly also
\[
\theta(x) = \sum_{\beta \in \bbZ^m} d_{\beta}e_{\beta}(x) \mbox{ with } \sum_{\beta \in \bbZ^m} |d_{\beta}| < \infty
\]
if $x \in [-L,L]^m$. If $(x,y) \in \left[ -\frac {L}{2},\frac {L}{2}\right]^{2m}$ we then have
\[
\begin{aligned}
\varphi(x-y)K(x-y)\theta(x) &= \sum_{(\alpha,\beta) \in \bbZ^{2m}} c_{\alpha}e_{\alpha}(x-y)d_{\beta}e_{\beta}(x) \\
&= \sum_{(\alpha,\beta) \in \bbZ^{2m}} c_{\alpha}d_{\beta}e_{\alpha+\beta}(x)e_{-\beta}(y) \\
&= \sum_{(\alpha,\beta) \in \bbZ^{2m}} c_{\alpha+\beta}d_{-\beta}e_{\alpha}(x)e_{\beta}(y)
\end{aligned}
\]
and
\[
\sum_{(\alpha,\beta) \in \bbZ^{2m}} |c_{\alpha+\beta}d_{-\beta}| = \sum_{(\alpha,\beta) \in \bbZ^{2m}} |c_{\alpha}|\,|d_{\beta}| < \infty.
\]
On the other hand, the operator on $L^2\left( \left[ -\frac {L}{2},\frac {L}{2}\right]^m,d\lambda\right)$ defined by the kernel which is the restriction of $e_{\alpha}(x)e_{\beta}(y)$ is rank~1 and its trace-norm is a constant $B > 0$ which depends only on $L$ and $m$. It follows that the trace-norm of the operator we are estimating and which is defined by a kernel with support in $\left[ -\frac {L}{2},\frac {L}{2}\right]^2$ is $\le \sum_{(\alpha,\beta) \in \bbZ^{2m}} |c_{\alpha}|\,|d_{\beta}| \cdot B < \infty$.\hfill\qed

\bigskip
\noindent
{\bf 9.2. Lemma.} {\em There is a constant $C > 0$ so that
\[
|S_{n,l}|_1 \le C \cdot 2^{(-1+1/p_1)\cdot n}.
\]
}

\bigskip
\noindent
{\bf {\em Proof.}} The operator $S_{n,l}$ is a compression of the operator defined by the $K_{n,l}$ in $L^2(\bbR^m,d\lambda)$. Denoting by ${\tilde S}_{n,l}$ the latter operator, it will suffice to show that $|{\tilde S}_{n,l}|_1$ satisfies the same kind of estimates. Note also that changing $l$ to $0$ corresponds to shifting in $\bbR^m$ by $\sigma(2^{-n})l$ and hence $|{\tilde S}_{n,l}|_1 = |{\tilde S}_{n,0}|_1$.

On the other hand, we have
\[
K_{n,0}(x,y) = K_{0,0}(\sigma(2^n)x,\sigma(2^n)y) \cdot 2^{n/p_1}
\]
and because $\int_{\bbR^m} f(\sigma(t)x)d\lambda(x) = t^{-1} \int_{\bbR^m} f(x)d\lambda(x)$ it is easy to see that
\[
|{\tilde S}_{n,0}|_1 = 2^{n/p_1} \cdot 2^{-n}|{\tilde S}_{0,0}|_1 = 2^{(-1+1/p_1)n}|{\tilde S}_{0,0}|_1
\]
and we may take $C = |{\tilde S}_{0,0}|_1$ which is finite by the preceding Lemma.\hfill\qed

\bigskip
We return to the proof of Proposition~9.2.

To satisfy the orthogonality conditions on the pieces into which we split our operators, we will need to replace $S_n$ by $\sum_{l \in g} S_{n,l}$ where $g \in \bbZ^m/(10^3\bbZ)^m \simeq (\bbZ/10^3\bbZ)^m$. If $l_1,l_2 \in g$ and $l_1 \ne l_2$, then $\|l_1-l_2\| \ge 10^3$ and $\theta(\sigma(2^n)x-10^{-2}l_1)\theta(\sigma(2^n)x-10^{-2}l_2) = 0$ since $\|10^{-2}l_1-10^{-1}l_2\| \ge 10$ and $\mbox{supp } \theta \subset B(0,10^{-1})$. Thus insures that $S^*_{n,l_1} S_{n,l_2} = 0$ since we have $K_{n,l_1}(x,y)K_{n,l_2}(x,z) = 0$. Since $\mbox{supp } \varphi \subset D(0,2) \subset B(0,2)$ we will have $K_{0,l}(x,y) = 0$ if $\|x-y\| \ge 2$. Thus if $K_{0,l_1}(x,y)K_{0,l_2}(y,z) \ne 0$ we must have $\|x-y\| \le 2$, $\|y-z\| \le 2$. On the other hand we must also have $\theta(x-l_1)\theta(y-l_2) \ne 0$ which is satisfied only if $\|x-l_1\| \le 2$, $\|y-l_2\| \le 2$ so that $\|l_1-l_2\| \le 6$ which is not satisfied if $l_1,l_2 \in g$, $l_1 \ne l_2$. Applying $\sigma(2^n)$ to $x,y,z$ we get $K_{n,l_1}(x,y)K_{n,l_2}(y,z) = 0$. We must still show that $K_{0,l_1}(x,y)K_{0,l_2}(z,y) = 0$. If this product is $\ne 0$ then $\|x-y\| \le 2$, $\|y-z\| \le 2$ and $\theta(x-l_1)\theta(z-l_2) \ne 0$ which implies that $\|x-l_1\| \le 2$, $\|z-l_2\| \le 2$ and hence $\|l_1-l_2\| \le 8$, which fails if $l_1,l_2 \in g$ and $l_1 \ne l_2$. Thus $K_{0,l_1}(x,y)K_{0,l_2}(z,y) = 0$. Applying $\sigma(2^n)$ to $x,y,z$ we get $K_{n,l_1}(x,y)K_{n,l_2}(z,y) = 0$.

Next we need to estimate how many of the operators $S_{n,l}$ where $l \in g$ are $\ne 0$ for a given $n$. Since the support of the measure $\mu$ is $[-1,1]^m$, in order that $S_{n,l} \ne 0$ we must have $(\mbox{supp } K_{n,l}) \cap [-1,1]^{2m} \ne \emptyset$. On the other hand $K_{n,l}(x,y) = K_{0,l}(\sigma(2^n)x,\sigma(2^n)y)$ so we must have $(\mbox{supp } K_{0,l}) \cap (\sigma(2^n)D(0,1))^2 \ne \emptyset$. Since $K_{0,l_1}(x,y)K_{0,l_2}(x,y) = 0$ the supports are disjoint and if $\lambda_{2m}(\mbox{supp } K_{0,l}) = h$ we infer that
\[
\begin{aligned}
&\mbox{card}\{l \in g \mid (\mbox{supp } K_{0,l}) \cap (\sigma(2^n)D(0,1))^2 \ne \emptyset\} \cdot h \\
&\quad = \mbox{card}\{l \in g \mid (\mbox{supp } K_{0,l}) \cap (D(0,2^n))^2 \ne \emptyset\} \cdot h \\
&\quad \le (\lambda_m(D(0,2^n)+B(0,R)))^2
\end{aligned}
\]
where $R$ is an upper bound on the diameter of $\mbox{supp } K_{0,l}$. On the other hand if $n \ge 0$ the right-hand side of the last inequality is $(2^{n/p_1+1}+R) \dots (2^{n/p_n+1}+R)2^m \le 2^{n+m(2+R)}$. It follows that
\[
\begin{aligned}
&\mbox{card}\{l \in g \mid \mbox{supp } K_{0,l} \cap (\sigma(2^n)D(0,1))^2 \ne \emptyset\} \\
&\quad \le h^{-1}2^{n+m(2+R)} \le C \cdot 2^n
\end{aligned}
\]
for some constant $C$.

Since $S_{n,l} = 0$ if $n \le n_0 - 2$, leaving out the $S_{n,l}$ with $n \le 0$ means leaving out only a finite number of trace-class operators. Thus we can apply he Lemma with $n \ge 0$, $\gamma_n = C2^n$ (to avoid $\gamma_n \le C2^n$ we may create a few zero terms). Let $q_1 = \frac {p_1}{p_1-1}$. Then if $\gamma_n = C_12^n$, $\beta_n = C_22^{-n/q_1}$ and $\pi_j  = j^{-1/q_1}$ we can choose $C$ large enough so that
\[
\gamma_{n+1} \le C\gamma_n,\ \sum_{k > n} \beta_k \le C\beta_n,\ \sum_{1\le k\le n} \beta_k\gamma_k \le C\beta_n\gamma_n.
\]
For the last inequality, the left-hand side is
\[
\le \mbox{const. } 2^{(n+1)(1-1/q_1)} \le \mbox{const. } 2^n \cdot 2^{-n/q_1}.
\]
Since $\beta_n = C_22^{-n/q_1}$, while $\pi_{\gamma_n} = (C_12^n)^{-n/q_1}$ we get $\pi_{\gamma_n} \ge C_3\beta_n$. This concludes the proof that the assumptions of Proposition~9.1 are satisfied, so we conclude that the proof of Proposition~9.2 by applying Proposition~9.1.\hfill\qed

\section{The formula for $k_{\varphi}$ of a commuting $n$-tuple of hermitian operators}
\label{sec10}

After having dealt with the technical problems posed by singular integrals with mixed homogeneities for the generalization to hybrid normed ideal perturbations, we can now copying the arguments in \cite{13} easily obtain the formula for $k_{\varphi}$ for commuting $m$-tuples of hermitian operators.

We begin with the immediate basic consequence of Proposition~9.2.

\bigskip
\noindent
{\bf 10.1. Proposition.} {\em Let $p_j > 1$, $1 \le j \le m$ be so that $\sum_{1 \le j \le m} p_j^{-1}=1$, and let $\varphi \in \cF([m])$ be given by $\varphi(j) = \Phi^-_{p_j}$, $1 \le j \le m$. Then if $\delta$ is the $m$-tuple of multiplication operators by the coordinate functions in $L^2([-1,1]^m,d\lambda)$ we have $k_{\varphi}(\delta) > 0$.
}

\bigskip
\noindent
{\bf {\em Proof.}} Let $X_j$ be the operator in $L^2([-1,1]^m,d\lambda)$ given by the kernel $K_j(x-y)$ where
\[
K_j(x_1,\dots,x_m) = (\mbox{sign } x_j)|x_j|^{p_j-1}(|x_1|^{p_1}+\dots + |x_m|^{p_m})^{-1}.
\]
Then by the technical results of Sections~8 and 9 we have $X_j \in \cC_{q_j}^+$ where $p_j^{-1} + q_j^{-1} = 1$. On the other hand, if $\delta = (D_1,\dots,D_m)$ is the $m$-tuple of multiplication operators by the coordinate functions in $L^2([-1,1]^m,d\lambda)$ we have that $\sum_{1 \le j \le m} [D_j,X_j]$ is the operator given by the kernel
\[
\sum_{1 \le j \le m} (x_j-y_j)K_j(x-y) = 1.
\]
That is $\sum_{1 \le j \le m} [D_j,X_j]$ is an operator of rank one proportional to the projection onto the one-dimensional subspace of constant functions in $L^2([-1,1]^m,d\lambda)$. By Proposition~4.1 this implies $k_{\varphi}(\delta) > 0$ since $\mbox{Tr } \sum_{1 \le j \le m} [D_j,X_j] > 0$.\hfill\qed

\bigskip
\noindent
{\bf 10.2. Proposition.} {\em Let $\varphi \in \cF([m])$ and let $\delta$ be the $m$-tuple of multiplication operators by the coordinate functions in $L^2([-1,1]^m,d\lambda)$. Then if $k_{\varphi}(\delta) < \infty$, there is a universal constant $\gamma_{\varphi}$ so that if $\tau$ is an $m$-tuple of commuting hermitian operators with Lebesgue absolutely continuous spectral measure, then we have
\[
(k_{\varphi}(\tau))^m = \gamma_{\varphi} \int m(x)d\lambda(x)
\]
where $m(x)$ is the multiplicity function of $\tau$. In particular, $\gamma_{\varphi} > 0$ iff $k_{\varphi}(\delta) > 0$.
}

\bigskip
\noindent
{\bf {\em Proof.}} The proof is in essence a repetition of the proof of Theorem~4.5 in \cite{13}. That is we first prove the formula when $\tau$ is the $m$-tuple of multiplication operators in $L^2(\Omega,d\lambda)$ where $\Omega$ is a bounded Borel set. Using Proposition~3.7 and Proposition~3.5 we prove this in stages: first for $\Omega$ a disjoint union of cubes, then for $\Omega$ an open set and then in full generality by using $K \subset \Omega \subset G$ where $K$ is compact and $G$ open with $\lambda(G\backslash K) < \epsilon$. We take into account that $k_{\varphi}(\tau)$ is invariant under shifts of $\Omega$ and homogeneous of degree~$1$ under dilations. The case of bounded multiplicity is then inferred from the result obtained, which is the case of multiplicity $\le 1$, by using Proposition~3.7. The general multiplicity case follows from the bounded multiplicity case using Proposition~3.5.\hfill\qed

\bigskip
\noindent
{\bf 10.1. Theorem.} {\em Let $p_j > 1$, $1 \le j \le m$ be so that $\sum_{1 \le j \le m} p_j^{-1} = 1$ and let $\varphi \in \cF([m])$ be defined by $\varphi(j) = \Phi_{p_j}^-$, $1 \le j \le m$. Then there is an universal constant $0 < \gamma < \infty$ depending only on the $p_j$, $1 \le j \le m$ so that if $\tau$ is an $m$-tuple of commuting hermitian operators then
\[
(k_{\varphi}(\tau))^m = \gamma \int m(x)d\lambda(x)
\]
where $m$ is the multiplicity function of the Lebesgue absolutely continuous part of the spectral measure of $\tau$.}

\bigskip
\noindent
{\bf {\em Proof.}} First we recall that by Proposition~7.4 $k_{\varphi}(\tau) = 0$ if the spectral measure of $\tau$ is singular with respect to Lebesgue measure. Then by Proposition~3.5 $k_{\varphi}(\tau)$ is equal to $k_{\varphi}(\tau_{ac})$ where $\tau_{ac}$ is the Lebesgue absolutely continuous part of $\tau$. This reduces the proof to the case where the spectral measure of $\tau$ is Lebesgue absolutely continuous. The only thing which is required in addition to Proposition~10.2 is that $0 < k_{\varphi}(\delta) < \infty$, where $\delta$ is the $m$-tuple of multiplication operators by the coordinate functions in $L^2([-1,1]^m,d\lambda)$. The $> 0$ fact is in Proposition~10.1, while $< \infty$ was shown in Proposition~7.4.\hfill\qed

\bigskip
The preceding theorem also easily gives that the generalized decomposition into $\varphi$-singular and $\varphi$-absolutely continuous subspaces in the case of $\varphi$ as in Theorem~10.1 and an $m$-tuple of commuting hermitian operators coincides with the decomposition into singular and absolutely continuous subspaces with respect to $m$-dimensional Lebesgue measure.

\section{Commutants\,$\mod \varphi$}
\label{sec11}

In this section we extend the definition of commutants $\mod$ normed ideals to the hybrid setting and extend the relation between $k_{\varphi}$ and approximate units for the compact ideal \cite{17}.

If $\tau \in \cL([n])$ is an $n$-tuple of hermitian operators and $\varphi \in \cF([n])$ we shall consider the $*$-subalgebras of $\cL$
\[
\begin{aligned}
\cE(\tau;\varphi) &= \{X \in \cL \mid [\tau(k),X] \in \cG_{\varphi(k)},\ k \in [n]\} \\
\cE^0(\tau;\varphi) &= \{X \in \cL \mid [\tau(k),X] \in \cG_{\varphi(k)}^{(0)},\ k \in [n]\} \\
\cK(\tau;\varphi) &= \cK \cap \cE(\tau;\varphi) \\
\cK^0(\tau;\varphi) &= \cK \cap \cE^0(\tau;\varphi)
\end{aligned}
\]
endowed with the norm
\[
\||X\|| = \|X\| + |[X,\tau]|_{\varphi}.
\]
If we drop the condition that $\tau = \tau^*$ we get non-self-adjoint operator algebras. There is also the following more general situation which is of interest: we may replace $\cF([n])$ by $\cF^*([n])$ where the $\varphi(k) \in \cF \cup \{\Phi_0\}$ with $\Phi_0$ denoting the ``norming function of the ideal $\{0\}$'' that is $|X|_{\Phi_0} \in \{0,\infty\}$ and $|X|_{\Phi_0} = 0 \Leftrightarrow X = 0$. In this case $\cE(\tau;\varphi)$ is equal $\cE(\tau;\varphi') \cap \{X_1,\dots,X_p\}'$, $\varphi' \in \cF([n])$ intersection of an $\cE(\tau;\varphi')$ with a von~Neumann algebra with a finitely generated commutant. Again, one may also consider the non-self-adjoint case.

\bigskip
\noindent
{\bf 11.1. Theorem.} {\em Let $\tau = \tau^* \in \cL([n])$ and let $\varphi \in \cF([n])$. Then we have
\begin{itemize}
\item[a)] the following are equivalent:
\begin{itemize}
\item[(i)] $k_{\varphi}(\tau) = 0$.
\item[(ii)] there are $A_m \in \cK^0(\tau;\varphi)$, $m \in \bbN$ so that
\[
\lim_{m \to \infty} \||X-A_mX\|| = \lim_{m \to \infty} \||X-XA_m\|| = 0
\]
for all $X \in \cK^0(\tau;\varphi)$ and $\||A_m\|| \le 1$, $m \in \bbN$.
\item[(iii)] condition {\rm (ii)} is satisfied and moreover $A_m \in \cR_1^+$ and $A_m \uparrow I$ as $m \to \infty$.
\end{itemize}
\item[b)] the following are equivalent:
\begin{itemize}
\item[(i)] $k_{\varphi}(\tau) < \infty$
\item[(ii)] there are $A_m \in \cK^0(\tau;\varphi)$, $m \in \bbN$ so that
\[
\lim_{m \to \infty} \||X-A_m\|| = \lim_{m \to \infty} \||X-XA_m\|| = 0
\]
for all $X \in \cK^0(\tau;\varphi)$ and $\sup_{m \in \bbN}\||A_m\|| < \infty$.
\item[(iii)] condition {\rm (ii)} is satisfied and moreover $A_m \in \cR_1^+$ and $A_m \uparrow I$ as $m \to \infty$.
\end{itemize}
\end{itemize}
}

\bigskip
\noindent
{\bf {\em Proof.}} a) (i) $\Rightarrow$ (iii) If $k_{\varphi}(\tau) = 0$ then there are $A_m \in \cR_1^+$, $A_m \uparrow I$ and $|[A_m,\tau]|_{\varphi} \to 0$ as $m \to \infty$ Clearly then $\|A_m\| \le 1$ and $\||A_m\|| \to 1$. Replacing the $A_m$'s by a subsequence and then passing from $A_m$ to $\alpha_mA_m$ for some suitable $\alpha_m \uparrow 1$ the new $A_m$'s will satisfy all the above conditions and additionally $\||A_m\|| < 1$.

If $X \in \cK^0(\tau;\varphi)$ then
\[
\lim_{m \to \infty} \|X-A_mX\| = \lim_{m \to \infty} \|X-XA_m\| = 0.
\]
Moreover, we have
\[
|[X-A_mX,\tau]|_{\varphi} \le |[A_m,\tau]|_{\varphi} \|X\|) + |(I-A_m)[X,\tau]|_{\varphi} \to 0
\]
as $m \to \infty$. Indeed $|[A_m,\tau]|_{\varphi} \to 0$ and $|(I-A_m)[X,\tau]_{\varphi} \to 0$ because $[X,\tau(k)] \in \cG_{\varphi(k)}^{(0)}$, $k \in [n]$. A similar argument gives
\[
|[X - XA_m, \tau]|_{\varphi} \to 0.
\]
Thus we have
\[
\lim_{m \to \infty} \||X-A_mX\|| = \lim_{m \to \infty} \||X-XA_m\|| = 0.
\]

Since (iii) $\Rightarrow$ (ii) trivially, we are left with (ii) $\Rightarrow$ (i).

Letting $X$ run over $\cR$ the fact that $\|X-XA_m\| + \|X-A_mX\| \to 0$ gives $A_m {\overset{s}{\to}} I$ and $A_m^* {\overset{s}{\to}} I$ as $m \to \infty$.  Hence $\liminf_{m \to \infty} \|A_m\| \ge 1$ so that $\||A_m\|| \le 1$ implies $\lim_{m \to \infty} |[A_m,\tau]|_{\varphi} = \lim_{m \to \infty} |[A_m^*,\tau]|_{\varphi} = 0$. This in turn gives $k_{\varphi}(\tau) = 0$ by an application of Proposition~3.4 to the sequence of hermitian operators $1/2(A_m+A_m^*) \in \cK^0(\tau;\varphi) \subset \cK$.

\bigskip
b) (i) $\Rightarrow$ (iii) If $k_{\varphi}(\tau) < \infty$ then there are $A_m \to \cR^+_1$, $A_m \uparrow I$ so that $\sup_{m \in \bbN} |[A_m,\tau]|_{\varphi} < \infty$. This then gives $\sup_{m \in \bbN} \||A_m\|| < \infty$. If $X \in \cK^0(\tau;\varphi)$ we then have
\[
\lim_{m \to \infty} \|X-XA_m\| = \lim_{m \to \infty} \|X-A_mX\| = 0.
\]
We also have
\[
|[(I-A_m)X,\tau]|_{\varphi} \le |(I-A_m)[X,\tau]|_{\varphi} + |[A_m,\tau]X_0|_{\varphi} + |[A_m,\tau]|_{\varphi}\|X-X_0\|
\]
where $X_0 \in \cR$. Since $[X,\tau(k)] \in \cG_{\varphi(k)}^{(0)}$ we have $\lim_{m \to 0} |(I-A_m)[X,\tau(k)]|_{\varphi(k)} = 0$ so that $|(I-A_m)[X,\tau]_{\varphi} \to 0$ as $m \to \infty$. Also $[A_m,\tau] {\overset{s}{\to}} 0$ so that $\|[A_m,\tau]X_0\| \to 0$ since $X_0 \in \cR$ and since 
\[
\mbox{rank}[A_m,\tau(k)]X_0 \le \mbox{rank } X_0
\]
we have $|[A_m,\tau]X_0|_{\varphi} \to 0$ as $m \to \infty$. Thus we have
\[
\limsup_{m \to \infty} |[(I-A_m)X,\tau]|_{\varphi} \le \left( \sup_{m \in \bbN} |[A_m,\tau]_{\varphi}\right) \cdot \|X-X_0\|.
\]
Since $X_0$ can be chosen so that $\|X-X_0\|$ be as close to zero as we want we get
\[
\lim_{m \to \infty} |[(I-A_m)X,\tau]|_{\varphi} = 0
\]
and hence $\||(I-A_m)X\|| \to 0$. Applying this to $X^*$ instead of $X$ gives also $\||X(I-A_m)\|| \to 0$.

Again (iii) $\Rightarrow$ (ii) is trivial and to prove that (ii) $\Rightarrow$ (i) like in the case of a) we infer from
\[
\|(I-A_m)X\| + \|X(I-A_m)\| \to 0
\]
that $A_m {\overset{s}{\to}} I$ and $A_m^* {\overset{s}{\to}} I$ as $m \to \infty$. Also from $\sup_{m \in \bbN} \||A_m\|| < \infty$ we get $\sup_{m \in \bbN}|[A_m,\tau]|_{\varphi} = \sup_{m \in \bbN} |[A_m^*,\tau]|_{\varphi} < \infty$. We then apply Proposition~3.4 to the sequence $1/2(A_m+A_m^*) \in \cK^0(\tau;\varphi) \subset \cK$ and we get that
\[
k_{\varphi}(\tau) < \infty.
\]
\hfill\qed

\bigskip
A similar result can be obtained for $qd_{\varphi}(\tau)$ working with finite rank projections instead of finite rank contractions.

\bigskip
\noindent
{\bf 11.2. Theorem.} {\em Let $\tau = \tau^* \in \cL([n])$ and $\varphi \in \cF([n])$. Then we have:
\begin{itemize}
\item[a)] The following conditions are equivalent:
\begin{itemize}
\item[(i)] $qd_{\varphi}(\tau) = 0$.
\item[(ii)] there are $P_m \in \cP$ so that
\[
\lim_{m \to \infty} \||(I-P_m)X\|| = 0
\]
for all $X \in \cK^0(\tau;\varphi)$ and
\[
\lim_{m \to \infty} \||P_m\|| = 1
\]
\item[(iii)] condition (ii) with the additional requirement $P_m \uparrow I$ as $m \to \infty$.
\end{itemize}
\item[b)] The following conditions are equivalent
\begin{itemize}
\item[(i)] $qd_{\varphi}(\tau) < \infty$.
\item[(ii)] there are $P_m \in \cP$ so that
\[
\lim_{m \to \infty} \||(I-P_m)X\|| = 0
\]
for all $X \in \cK^0(\tau;\varphi)$ and
\[
\sup_{m \in \bbN} \||P_m\|| < \infty.
\]
\item[(iii)] condition (ii) with the additional requirement that $P_m \uparrow I$ as $m \to \infty$.
\end{itemize}
\end{itemize}
}

\bigskip
We leave the proof, which is similar to that of Theorem~1.11 as an exercise for the reader.

\end{document}